\documentclass[pdflatex,sn-mathphys-num]{sn-jnl}

\usepackage{amsmath,amssymb,amsfonts}%
\usepackage{amsthm}%
\usepackage{mathrsfs}%
\usepackage[title]{appendix}%
\usepackage{xcolor}%
\usepackage{cases}%
\usepackage{textcomp}%
\usepackage{manyfoot}%
\usepackage{booktabs}%
\usepackage{algorithm}%
\usepackage{algorithmicx}%
\usepackage{algpseudocode}%
\usepackage{listings}%
\usepackage{lmodern}%
\usepackage{anyfontsize}%
\usepackage{lipsum}
\usepackage{graphicx}%
\usepackage{multirow}%
\usepackage{tikz}
\usepackage{float}
\usetikzlibrary{arrows.meta}
\usetikzlibrary{positioning}
\hypersetup{hypertexnames=false}

\newtheorem{theorem}{Theorem}[section]
\newtheorem{lemma}[theorem]{Lemma}%

\allowdisplaybreaks
\begin{document}

\title[A Littlewood-Paley approach to the Mittag-Leffler function\dots]{\large A Littlewood-Paley approach to the Mittag-Leffler function in the frequency space and applications to nonlocal problems}


\author{\fnm{\normalsize Ahmed A.} \sur{\normalsize Abdelhakim}}\email{\normalsize a.mouhamed@qu.edu.sa}\email{\normalsize ahmed.abdelhakim@aun.edu.eg}
%
%
\affil{\orgdiv{\normalsize Department of Mathematics}, \orgname{\normalsize College of Science, Qassim University}, \state{\normalsize Buraydah}, \postcode{\normalsize 51452}, \country{\normalsize Saudi Arabia}}
\affil{\orgdiv{\normalsize Department of Mathematics}, \orgname{\normalsize Faculty of Science, Assiut University}, \state{\normalsize Assiut}, \postcode{\normalsize 71516}, \country{\normalsize Egypt}}
%
%


\abstract{Let $0<\alpha<2$, $\beta>0$ and $\alpha/2<|s|\leq 1$. We obtain values of the Lebesgue exponent $p=p(\gamma)$, $\gamma>0$, for which the Fourier transform of $ E_{\alpha,\beta}(
e^{\dot{\imath}\pi s} |\cdot|^{\gamma})$
is an $L^{p}(\mathbb{R}^d)$ function, using tools from the Littlewood-Paley theory. This question arises in the analysis of certain space-time fractional diffusion and Schr\"{o}dinger
problems and has been solved for the particular cases $\alpha\in (0,1)$, $\beta=\alpha,1$, and
$s=-1/2,1$ via asymptotic analysis of Fox $H$-functions. The Littlewood-Paley theory
provides a simpler proof that allows considering all values of $\beta,\gamma>0$ and $s\in (-1,1]\setminus [-\alpha/2,\alpha/2]$. This enabled us to prove various key estimates for a general class of nonlocal space-time problems.}




\maketitle

\section{Introduction}
The Mittag-Leffler function $E_{\alpha,\beta}$
is entire of order $1/\alpha$, when $\alpha,\beta>0$ (see e.g. \cite{gorenflo2020mittag}, Section 4.1). It is also bounded and decays like $|z|^{-1}$
in the sector $\{z\in \mathbb{C}: |\arg{z}|>\alpha \pi/2 \}$ (see e.g. \cite{podlubny1998fractional}, Theorem 1.6). So, for all $s\in (-1,1]\setminus [-\alpha/2,\alpha/2]$ and all $\gamma>0$, its composition with the function $\mathbb{R}^{d}\ni x\mapsto e^{\dot{\imath} \pi s } |x|^{\gamma}$ is a continuous bounded function that decays like $|x|^{-\gamma}$. Consequently,
$E_{\alpha,\beta}\left(
e^{\dot{\imath}\pi s} |\cdot|^{\gamma}\right)$ has a radially symmetric Fourier transform in the sense of tempered distributions, for all $s\in (-1,1]\setminus [-\alpha/2,\alpha/2]$ and all $\gamma>0$. In fact, it has a Fourier transform as an $L^{p}(\mathbb{R}^{d})$ function when $d/\gamma<p\leq 2$, by the Hausdorff-Young inequality. The question of determining the integrability properties of the Fourier transform of
$E_{\alpha,\beta}\left(
e^{\dot{\imath}\pi s} |\cdot|^{\gamma}\right)$ is therefore more interesting when
$0<\gamma\leq d/2$.
\medskip \par
Let $\mathcal{F}f$  denote the Fourier transform of a function $f$. The following $L^{p}(\mathbb{R}^{d})$ properties of $\mathcal{F}
E_{\alpha,\beta}(e^{\dot{\imath}
s } |\cdot|^{\gamma})$ can be found in \cite{abdelhakim2025asymptotic}.
\begin{theorem}(\cite{abdelhakim2025asymptotic}, Theorem 3)\label{thpd}
Suppose that $\alpha\in (0,2)$ and $\beta>0$.
Fix $\,-\pi<s\leq \pi$ such that
$|s|>\alpha \pi/2$. Then
$\mathcal{F} E_{\alpha,\beta}\left(
e^{\dot{\imath}\pi s} |\cdot|^{\gamma}\right)\in
L^{p}(\mathbb{R}^{d})$, for all
\begin{equation*}
 \left\{
   \begin{array}{ll}
 1<p<{d}/{(d-\gamma)}, &\hspace{0.5cm} \hbox{if\hspace{0.1cm} $(d-1)/2<\gamma<d$;}
\medskip \\
    1<p<\infty, & \hspace{0.5cm} \hbox{if\hspace{0.1cm} $\gamma=d$;} \\[0.1 cm]
        1<p\leq\infty, & \hspace{0.5cm} \hbox{if\hspace{0.1cm} $\gamma>d$.}
   \end{array}
 \right.
\end{equation*}
\end{theorem}
Theorem \ref{thpd} is proved by showing that $\xi\mapsto\mathcal{F}E_{\alpha,\beta}(e^{\dot{\imath}\pi s} |\cdot|^{\gamma})(\xi)$ is continuous on $\mathbb{R}^d\setminus \{0\}$, and obtaining asymptotic upper bounds for it,
both when $\xi\to 0$ and when $\xi\to\infty$, provided that $\gamma>(d-1)/2$ (see \cite{abdelhakim2025asymptotic}, Theorem 2).
\smallskip\par  For smaller values of $\gamma$, $E_{\alpha,\beta}\left(
e^{\dot{\imath}\pi s}|\cdot|^{\gamma} \right)$ is less smooth at the origin and its Fourier transform expectedly shows slower decay at infinity. In this work, we use the Littlewood-Paley theory
to prove that the Fourier transform of
$E_{\alpha,\beta}\left(
e^{\dot{\imath}\pi s} |\cdot|^{\gamma}\right)$ is in $L^{p}(\mathbb{R}^{d})$ for certain $p=p(\gamma)$,
for every $\gamma>0$, including all the values of $p$ already obtained in Theorem \ref{thpd}. We thus recover the low regularity case $0<\gamma \leq (d-1)/2$ missing in Theorem \ref{thpd}. \smallskip
\begin{theorem}\label{lmlwp}
Let $\alpha \in (0,2)$ and $\beta >0$.
Suppose that $\gamma>0$ and $\alpha/2< |s|\leq 1$. Then
$\mathcal{F}E_{\alpha,\beta}(e^{\dot{\imath}\pi s}
|\cdot|^{\gamma} )
\in L^{p}(\mathbb{R}^{d})$
whenever $p>1$ and $\gamma p^{\prime}>d$. In other words, $\mathcal{F}
E_{\alpha,\beta}(e^{\dot{\imath}\pi s}
|\cdot|^{\gamma} )\in L^{p}(\mathbb{R}^{d})$ for all
\begin{equation*}
 \left\{
   \begin{array}{ll}
 1<p<{d}/{(d-\gamma)}, &\hspace{0.5cm} \hbox{if\hspace{0.1cm} $0<\gamma<d$;}
\medskip \\
    1<p<\infty, & \hspace{0.5cm} \hbox{if\hspace{0.1cm} $\gamma=d$;} \\[0.1 cm]
        1<p\leq\infty, & \hspace{0.5cm} \hbox{if\hspace{0.1cm} $\gamma>d$.}
   \end{array}
 \right.
\end{equation*}
\end{theorem}
The interest in these $L^p$ properties is due to the fact that solutions
to certain fractional heat and Schr\"{o}dinger problems can be represented as a convolution operator with the Fourier transform
of $E_{\alpha,\beta}(e^{\dot{\imath}\pi s}
|\cdot|^{\gamma})$ in some direction $s$.
We make this precise and give more details in Section \ref{apls}. Solutions to this type of equations can also be obtained using the Mellin transform in terms of
a convolution with Fox $H$-functions (see e.g. \cite{mainardi2005fox,kemppainen2017representation,su2019local}
and the monograph \cite{mathai2009h}).
This is the method adopted e.g. in \cite{kemppainen2017representation,su2019local}
where the authors show that the resulting
$H$-functions are bounded in  $L^{p}(\mathbb{R}^{d})$ via an asymptotic analysis. These $L^{p}(\mathbb{R}^{d})$ results obtained in \cite{kemppainen2017representation,su2019local} can be shown to be related to the results of Theorem \ref{lmlwp} when $\alpha\in (0,1)$ and $\beta=\alpha,1$, in the directions $s=1$ and $s=-1/2$, respectively.\medskip \par Let ${f}^{\vee}$
denote the inverse Fourier transform of a tempered distribution $f$. In order to prove Theorem \ref{lmlwp}, we apply the Littlewood-Paley projection operator $P_{j}$, $j\in \mathbb{Z}$, to the inverse Fourier transform of $E_{\alpha,\beta}(
e^{\dot{\imath}\pi s} |\cdot|^{\gamma})$. By the Littlewood-Paley theory, it suffices to show that the sum
\begin{equation*}
\sum_{j\in \mathbb{Z}} \| P_{j}\left(E_{\alpha,\beta}\left(
e^{\dot{\imath}\pi s} |\cdot|^{\gamma}\right)\right)^{\vee}\,
\|_{L^{p}(\mathbb{R}^{d})}
\end{equation*}
is convergent.
The main task then is to estimate the $L^{p}(\mathbb{R}^{d})$ norm of these projections effectively. Our starting point is the formula
\begin{equation}\label{int1}
P_{j}\left(E_{\alpha,\beta}(
e^{\dot{\imath}\pi s}
 |\cdot|^{\gamma} )\right)^{\vee}(\xi)
=\frac{2\pi}{|\xi|^{\frac{d}{2}-1}}
\int_{0}^{\infty}\psi(2^{-j}r)
E_{\alpha,\beta}(
e^{\dot{\imath}\pi s} r^{\gamma} ) J_{\frac{d}{2}-1}(2\pi |\xi| r)
r^{\frac{d}{2}}dr,
\end{equation}
where $J_{{d}/{2}-1}$ is the Bessel function and $\psi$ is a smooth positive bump function supported in $[1/2,2]$. Upon substituting for $J_{d/2-1}$ into (\ref{int1}), each projection can be written down as a radial oscillatory integral. Roughly speaking, we break up this oscillatory integral to separate low and high frequency components of the projection. The latter appear as highly oscillating slowly decaying integrals, and are therefore more difficult to estimate. The proof of Theorem \ref{lmlwp} is detailed in Section \ref{scprf}.\par
In Section \ref{apls}, we apply Theorem \ref{lmlwp} to a general class of nonlocal space-time problems that include the space-time fractional heat and Schr\"{o}dinger equations. We prove dispersive estimates and decay estimates for the $L^{p}(\mathbb{R}^{d})$ norm of the solution in its spacial variable for that class. These estimates are key to
studying the convergence of the solution to the initial data as well as to questions of well-posedness of the corresponding semilinear equations (see e.g. \cite{kemppainen2017representation, li2022applicationi,li2022applicationii,su2019local}).

\section{Preliminaries}
\subsection{ A contour integral representation for the derivatives of the Mittag-Leffler function}
Let $\alpha\in (0,2)$, $\beta\in \mathbb{C}$
and $\rho>0$. Fix $\pi \alpha/2<\omega< \min\{\pi \alpha,\pi\}$ and let $C_{\rho, \omega}$
be the contour that consists of the two rays
$\left\{z\in \mathbb{C}:
\arg{z}=\pm\omega, |z|\geq \rho\right\}$ along with the circular arc $\left\{z\in \mathbb{C}:-\omega\leq\arg{z}\leq \omega, |z|=\rho\right\}$, taken as in Fig. 1 below. Then, the contour integral
\begin{equation}\label{cntr0}
\frac{1}{2\pi \dot{\imath} \alpha}\int_{C_{\rho,\omega}}
\frac{e^{\zeta^{1/\alpha}}{
\zeta^{(1-\beta)/\alpha}}}{
\zeta-z}d\zeta
\end{equation}
is analytic and converges in the sector $|\arg{z}|>\omega$ to the Mittag-Leffler function \begin{equation*}
E_{\alpha,\beta}(z)=\sum_{k=0}^{\infty}
\frac{z^{k}}{\Gamma{(\alpha k+ \beta)}}.
\end{equation*}
The representation (\ref{cntr0}) is derived originally in \cite{dzhrbashyan1954integral}.
See also \cite{abdelhakim2025asymptotic}, Section 2 and \cite{gorenflo2020mittag}, Section 4.7.
Given $r\geq 0$ and
$-1<s\leq 1$ such that
$\,|s|> \alpha/2$, we can set
\begin{equation*}
\alpha \pi/2<\omega<\min
\left\{|s|,\alpha \right\}\pi,\; \rho=1,
\end{equation*}
and use the contour integral (\ref{cntr0}) to obtain
\begin{equation}\label{esln1}
E_{\alpha,\beta}(r e^{\dot{\imath} \pi s })=
\frac{1}{2\pi \dot{\imath} \alpha}\int_{C_{\omega}}
\frac{e^{\zeta^{1/\alpha}}{
\zeta^{(1-\beta)/\alpha}}}{
\zeta-r e^{\dot{\imath} \pi s }}d\zeta,
\end{equation}
where, for short, $C_{\omega}$ denotes the contour $C_{1,\omega}$. See Figure \ref{fig1} below. By analyticity, it follows from
(\ref{esln1}) that the $k^{\text{th}}$ derivative
of $E_{\alpha,\beta}(r e^{\dot{\imath} \pi s})$ equals
\begin{equation}\label{dme0}
\frac{k! e^{\dot{\imath} \pi k s}}{2\pi \dot{\imath} \alpha}\int_{C_{\omega}}
\frac{e^{\zeta^{1/\alpha}}{
\zeta^{(1-\beta)/\alpha}}}{
(\zeta-r e^{\dot{\imath} \pi s })^{k+1}}d\zeta.
\end{equation}
We are going to make use of the following lower bound
on the distance from a point on the ray $\arg{z}=\pi s$, $-1<s\leq 1$,
to the contour $C_{\omega}$ (cf. the inequality
(2.4) in \cite{abdelhakim2025asymptotic} and the inequality (4) in
\cite{abdelhakim2025mittag}):
\begin{equation}\label{nf1}
\inf_{\;\zeta\in C_{\omega}}{|\zeta-\eta e^{\dot{\imath} \pi s }|}\geq
|\sin{(\pi |s|-\omega)}|,\quad \eta\geq 0.
\end{equation}
This follows from the facts that
\begin{eqnarray*}
\inf_{\;\zeta\in C_{\omega}}{|\zeta-\eta e^{\dot{\imath} \pi s }|}
&=&
\inf_{\rho\geq 1,\,|\arg{\zeta}|\leq
 \omega}{|\rho e^{\dot{\imath}\arg{\zeta} }-\eta e^{\dot{\imath} \pi s }|}\\
&\geq&  \inf_{\rho\geq 1}{
\sqrt{\rho^2+\eta^2-2\rho\eta \cos{(\pi |s|-\omega)}}},
\end{eqnarray*}
and that
\begin{equation*}
\min_{\eta\geq 0}{\sqrt{\rho^2+\eta^2-2\rho\eta \cos{(\pi |s|-\omega)}}}=
\rho|\sin{(\pi |s|-\omega)}|.
\end{equation*}
Figure \ref{fig1} below gives an illustration of the inequality (\ref{nf1}).
\begin{figure}[H]
\centering
\begin{tikzpicture}[scale=0.8]
\draw [dotted,->](-4.5,0) --(4.5,0);
\draw [dotted,->](0,-4.5) --(0,4.5);
\draw [very thick,brown] (0.32,0.383) arc [radius=0.55, start angle=50, end angle= 110];
\node [very thick, brown, left]
 at (-0.3,0.3) {$\textstyle\theta=\pi |s|-\omega$};
\node [very thick, brown, above]
 at (0,0) {$\scriptstyle\mathbf{\theta}$};
 \node [very thick, rotate=-15,brown]
 at (0,2.7) {$\textstyle d(\eta e^{\dot{\imath} \pi|s|},C_{\omega})$};
\draw [thick,brown,dashed](-1.268,2.719)--(1.607,1.915);
\draw [thick,dotted,blue](0,0)--(1.607,1.915);
\draw [thick,blue](1.607,1.915)--(3.214,3.83);
\draw [-{Latex[length=3mm]}][blue](1.607,1.915)--(2.41,2.873);
\draw [thick,dotted,blue](0,0)--(1.607,-1.915);
\draw [thick][blue](1.607,-1.915)--(3.214,-3.83);
\draw [-{Latex[length=3mm]}][blue](3, -3.577)--(2.41,-2.873);
\draw [thick,blue] (2.5,0) arc [radius=2.5, start angle=0, end angle= 50];
\draw [thick,blue] (2.5,0) arc [radius=2.5, start angle=0, end angle= -50];
\draw [blue] [-{Latex[length=3mm]}] (2.266,1.057) arc [radius=2.5, start angle=25, end angle= 35];
\node [below right] at (2.4,0.1) {$\textstyle 1$};
\node [right,blue] at (2,1.5) {$\textstyle C_{\omega}$};
\draw [dashed][red](0,0)--(3.83,3.214);
\draw [dashed][red](0,0)--(3.83,-3.214);
\draw [thick,blue] (0.45,0) arc [radius=0.45, start angle=0, end angle= 50];
\node [above left] at (0.55,-0.1) {$\scriptstyle\textcolor{blue}{\omega}$};
\draw [thick,red] (0.40,0) arc [radius=0.40, start angle=0, end angle= -40];
\node at (0.9,-0.22) {$\textstyle\textcolor{red}{\alpha \pi/2}$};
\draw [thick](0,0)--(-1.902,4.078);
\draw [thick](0,0)--(-1.902,-4.078);
\draw [fill=black](-1.268,2.719) circle (.04 cm);
\node [left] at (-1.2,-2.6) {$\textstyle\arg{z}=-\pi|s|$};
\node [left] at (-0.634,1.36) {$\textstyle\eta$};
\end{tikzpicture}
\caption{The contour $C_{\omega}$ in blue is set up such that
$\pi \alpha/2<\omega<\min{\{|s|,\alpha\}}\pi$. The distance from
any point on the ray $\arg{z}=\pi s$, $-1<s\leq 1$,  to the contour $C_{\omega}$
is bounded below by $|\sin{(\pi |s|-\omega)}|$.}\label{fig1}
\end{figure}
\subsection{Asymptotics for Bessel functions}\label{bsl}
The Fourier transform of a radial function $\varphi:\mathbb{R}^{d}\rightarrow \mathbb{C}$
such that $\varphi(x)=\varphi_{0}(|x|)$
is the radial function given by
\begin{equation}\label{frtrdl}
\widehat{\varphi}(\xi)=
\frac{2\pi}{|\xi|^{\frac{d}{2}-1}}
\int_{0}^{\infty}
\varphi_{0}(r)J_{\frac{d}{2}-1}(2\pi |\xi| r)
r^{\frac{d}{2}}dr,
\end{equation}
where $J_{\lambda}$ is the Bessel function defined by
\begin{equation*}
J_{\lambda}(r):=
\frac{2^{-\lambda}}{\Gamma{(\frac{1}{2})}
\Gamma{(\lambda+\frac{1}{2})}}
r^{\lambda}
\int_{-1}^{1}e^{\dot{\imath}  r s}
(1-s^{2})^{\lambda-\frac{1}{2}}ds,\quad r\geq 0.
\end{equation*}
The identity (\ref{frtrdl}) is a consequence of radial symmetry in addition to the fact that ${2\pi}{|\xi|^{1-\frac{d}{2}}}
J_{{d}/{2}-1}(2\pi |\xi|)$ is the Fourier transform of the surface measure on the sphere
$\mathbb{S}^{d-1}$ (
see e.g. \cite{stein1993harmonic}, Chapter IX, and \cite{samko1993fractional}, Lemma 25.1). A sufficient condition for (\ref{frtrdl}) to hold is that
\begin{equation*}
\int_{0}^{\infty}r^{d-1}|\phi_{0}(r)|dr<\infty.
\end{equation*}
The function $J_{-1/2}$ can be computed explicitly. It is given by the identity
\begin{equation}\label{jd1}
J_{-\frac{1}{2}}(r)=
\sqrt{\frac{2}{\pi }}{r^{-\frac{1}{2}}}{\cos{r}},\quad r> 0.
\end{equation}
If $\operatorname{Re}\lambda>-1/2$, we have
\begin{equation*}
J_{\lambda}(r)=\frac{2^{-\lambda}}{
\Gamma{(\lambda+1)}} r^{\lambda}+O(r^{1+\operatorname{Re}
\lambda}),
\end{equation*}
which is obviously useful when the argument of
$J_{\lambda}$ is small. For $d>1$ in particular one has
\begin{equation}\label{jdsmll1}
J_{\frac{d}{2}-1}(r)=
\frac{2^{1-\frac{d}{2}}}{\Gamma{(\frac{d}{2})}}
r^{\frac{d}{2}-1}
+O(r^{\frac{d}{2}}).
\end{equation}
For more on Bessel functions, see e.g. \cite{grafakos2008classical}, Appendix B, and \cite{stein1993harmonic}, Chapter VIII. For a large argument, we have the asymptotic expansion:
Let $\operatorname{Re}\lambda>-1/2$. Then, for $r>1$, the Bessel function $J_{\lambda}(r)$
has the asymptotic expansion:
\begin{lemma}(\cite{abdelhakim2025asymptotic}, Lemma 3.1)
Suppose that $\operatorname{Re}\lambda>-1/2$. Then the Bessel function $J_{\lambda}(r)$
has the asymptotic expansion:
\begin{equation}\label{bssl}
J_{\lambda}(r)=
\sum_{0\leq \ell \leq M}
c_{\ell}(\lambda)\cos{(r+\lambda_{\ell})}
r^{-(\ell+\frac{1}{2})}+L_{\lambda}(r;M),
\end{equation}
for each $M\geq 1$, where
\begin{align*}
c_{\ell}(\lambda)&:=\,
\frac{\left(\lambda-\frac{1}{2}\right)_{\!\ell}}{\sqrt{2\pi} 2^{\ell-1}\ell !}
\frac{\Gamma{(\lambda+\ell+\frac{1}{2})}}{
\Gamma{(\lambda+\frac{1}{2})}},\\[0.1 cm]
\lambda_{\ell}&:=\,\frac{\pi}{2}(\ell-\lambda)-\frac{\pi}{4},
\end{align*}
and $r\mapsto L_{\lambda}(r;M)$ is a continuous function on $(1,\infty)$ such that
\begin{equation*}
\left|L_{\lambda}(r;M)\right|\lesssim
r^{-M-\frac{3}{2}}.
\end{equation*}
\end{lemma}
See \cite{abdelhakim2025asymptotic}, Appendix A, for a proof of (\ref{bssl}) where $L_{\lambda}(r;M)$ is computed explicitly.
\subsection{The Littlewood-Paley projection operator}
\par Fix $j\in \mathbb{Z}$ and a real-valued radially symmetric bump function $\varphi$ supported in the ball $\{x\in \mathbb{R}^{d}:|x|\leq 2\}$ and equals $1$ on the unit ball. The Littlewood-Paley projection operator  $P_{j}$ is the Fourier multiplier defined by
\begin{equation*}
\widehat{P_{j}f}(\xi):=(\varphi
({\xi}/{2^{j}})-
\varphi({\xi}/{2^{j-1}}))
\widehat{f}(\xi).
\end{equation*}
The operator $P_{j}$ is by definition a smoothed out projection to the lacunary set $\{x\in \mathbb{R}^{d}:2^{j-1}\leq |x|\leq 2^{j+1}\}$.
The well-known Littlewood-Paley inequality reads
\begin{equation*}
\|f\|_{L^{p}(\mathbb{R}^{d})}\approx
\|(\sum_{j\in \mathbb{Z}} |P_{j}f|^{2})^{1/2}\|_{L^{p}(\mathbb{R}^{d})},\quad
1<p<\infty.
\end{equation*}
For the basic Littlewood-Paley theory and its wide applications see e.g. \cite{grafakos2008classical,stein1993harmonic,tao2006nonlinear}.
We shall make use of the following inequality which turns out to be sufficient for our purposes:
\begin{equation}\label{plt}
\|f\|_{L^{p}(\mathbb{R}^{d})}\lesssim
\sum_{j\in \mathbb{Z}} \| P_{j}f\|_{L^{p}(\mathbb{R}^{d})},\quad 1< p<\infty.
\end{equation}
\section{Proof of Theorem \ref{lmlwp}}\label{scprf}
Fix $j\in \mathbb{Z}$. Using (\ref{int1}) and changing variables $2^{-j+1}\pi |\xi|r
\rightarrow r$, we see that
\begin{equation}\label{plttrngl}
P_{j}\left(E_{\alpha,\beta}(
e^{\dot{\imath}\pi s}
 |\cdot|^{\gamma} )\right)^{\vee}(\xi)
=(2\pi)^{-\frac{d}{2}}
\left({\mathcal{L}}_{j}(\xi)+
{\mathcal{H}}_{j}(\xi)\right),
\end{equation}
where
\begin{eqnarray*}
{\mathcal{L}}_{j}(\xi) &:=& \int_{\left\{0\leq r\leq 2^{-j}\right\}}Q_{j}(r,\xi)dr, \\
{\mathcal{H}}_{j}(\xi)&:=& \int_{\left\{r> 2^{-j}\right\}}Q_{j}(r,\xi)dr,
\end{eqnarray*}
with
\begin{eqnarray*}
Q_{j}(r,\xi)&:=&\frac{2^{(\frac{d}{2}+1)j}}
{|\xi|^{d}}
\Phi_{j}(r,\xi)
J_{\frac{d}{2}-1}
\left({2^{j}r}\right)
r^{\frac{d}{2}},\\
\Phi_{j}(r,\xi)&:=&
  \psi\left(\frac{r}{2\pi |\xi|}\right)
E_{\alpha,\beta}\left( e^{\dot{\imath}\pi s} \left(\frac{2^{j}r}{2\pi  |\xi|}\right)^{\gamma}\right),
\end{eqnarray*}
and $\psi$ is a positive test function with support in $[1/2,2]$. \par
The splitting (\ref{plttrngl})
is natural considering the difference between the asymptotic behaviour of the Bessel function $J_{{d}/{2}-1}$ of a small argument
and its asymptotic behaviour when its argument is large, as demonstrated by (\ref{jd1}) for $d=1$, and by (\ref{jdsmll1}) and (\ref{bssl}), for $d>1$. Most importantly, this decomposes the projection $P_{j}\left(E_{\alpha,\beta}(
e^{\dot{\imath}\pi s} |\cdot|^{\gamma} )\right)^{\vee}$ into the sum of its low and high frequency components.
Notice that, up to some universal constant, $\mathcal{L}_{j}$ has frequency smaller
than $2^{-j}$ while $\mathcal{H}_{j}$ has frequency larger than $2^{-j}$. The Mittag-Leffler function $E_{\alpha,\beta}$ has the property that,  whenever $|s|>\alpha/2$, $E_{\alpha,\beta}(
e^{\dot{\imath}\pi s} 2^{\gamma j})$
is uniformly bounded for all $j\in \mathbb{Z}$ and decays like $2^{-\gamma j}$ for large $j>0$ (see e.g. \cite{abdelhakim2023asymptotic} and the references therein). The plan is to estimate the $L^{p}(\mathbb{R}^{d})$ norms of
both the low frequency component ${\mathcal{L}}_{j}$ and the high frequency component
${\mathcal{H}}_{j}$, and then use the triangle inequality in the light of (\ref{plttrngl}) to obtain an estimate for $\| P_{j}\left(E_{\alpha,\beta}(
e^{\dot{\imath}\pi s} |\cdot|^{\gamma} )\right)^{\vee}
\|_{L^{p}(\mathbb{R}^{d})}$. $L^{p}$ boundedness of $\left(E_{\alpha,\beta}(e^{\dot{\imath}\pi s}
|\cdot|^{\gamma} )\right)^{\vee}$ follows then
by the Littlewood-Paley inequality (\ref{plt}). Before we proceed, we introduce the number
\begin{equation*}
\delta(j):=\begin{cases}
          -\gamma j, & j>0;\\[0.1 cm]
          0, & j\leq 0.
        \end{cases}
\end{equation*}
\par Let us start with ${\mathcal{L}}_{j}$. Suppose that $ 0\leq r\leq 2^{-j}$. Then,
by formula (\ref{jd1})  for $d=1$
and (\ref{jdsmll1}) for $d>1$, we have the estimate
\begin{equation}\label{jssl}
\left|J_{\frac{d}{2}-1}\left({2^{j}r}\right)\right|
\lesssim
{2^{(\frac{d}{2}-1)j}r^{\frac{d}{2}-1}}.
\end{equation}
Furthermore, we have that
\begin{equation}\label{stmml1}
\left|\Phi_{j}(r,\xi)\right|\lesssim
{2^{\delta(j)}}\,
\psi\left(\frac{r}{2\pi |\xi|}\right),
\end{equation}
which follows from the decay estimate (see e.g. \cite{abdelhakim2023asymptotic} and \cite{kilbas2006theory}, Theorem 4.3):
\begin{equation*}
|E_{\alpha,\beta}( r e^{\dot{\imath}\pi s})|\lesssim \frac{1}{1+r}, \quad r\geq 0,\;|s|>\alpha/2.
\end{equation*}
Combining (\ref{jssl}) and (\ref{stmml1}) we deduce that
\begin{eqnarray}
\nonumber |{\mathcal{L}}_{j}(\xi)|&\lesssim&\frac{2^{j d}}
{|\xi|^{d}}
\int_{\left\{0\leq r\leq 2^{-j}\right\}}
\left|\Phi_{j}(r,\xi)\right| r^{d-1}dr\\
\label{ljx0}&\lesssim&\frac{2^{j d+\delta(j)}}
{|\xi|^{d}}
\int_{\left\{0\leq r\leq 2^{-j}\right\}}\psi\left(\frac{r}{2\pi |\xi|}\right) r^{d-1}dr.
\end{eqnarray}
Invoking Minkowski's integral inequality and the fact that $\psi$ is supported in $[1/2,2]$, we obtain that
\begin{eqnarray*} \|{\mathcal{L}}_{j}\|_{L^{p}(\mathbb{R}^{d})}
&\lesssim&
{2^{j d+\delta(j)}} \int_{\left\{0\leq r\leq 2^{-j}\right\}}\left(\int_{\mathbb{R}^{d}}
{|\xi|^{-d p}}\psi\left(\frac{r}{2\pi|\xi|}\right)^{p}\,d\xi
\right)^{\frac{1}{p}}
r^{d-1}dr\\
&\leq& {2^{j d+\delta(j)}} \int_{\left\{0\leq r\leq 2^{-j}\right\}}\left(\int_{
\left\{r/4\leq \pi |\xi|\leq r\right\}}
{|\xi|^{-d p}}\,d\xi
\right)^{\frac{1}{p}}
r^{d-1}dr\\
&\lesssim& {2^{j d+\delta(j)}} \int_{\left\{0\leq r\leq 2^{-j}\right\}}
r^{\frac{d}{p}-1}dr\\
&\lesssim& {2^{d(1-\frac{1}{p})j+\delta(j)}},
\end{eqnarray*}
for all $1\leq p<\infty$. This estimate extends to $p=\infty$. Indeed, it follows from (\ref{ljx0}) that
\begin{equation*}
|{\mathcal{L}}_{j}(\xi)|\lesssim
\frac{2^{j d+\delta(j)}}
{|\xi|^{d}}
\int_{\left\{0\leq r\leq 2^{-j}\right\}\cap \left\{\pi |\xi|\leq r\leq 4\pi |\xi|\right\}} r^{d-1}dr\,\lesssim\,{2^{j d+\delta(j)}}.
\end{equation*}
We have so far shown that
\begin{equation}\label{nrmstl1}
\|{\mathcal{L}}_{j}\|_{L^{p}(\mathbb{R}^{d})}
\lesssim {2^{d(1-\frac{1}{p})j+\delta(j)}},
\end{equation}
for all $1\leq p\leq \infty$.\par
We turn our attention to the high frequency component ${\mathcal{H}}_{j}$. For the dimensions $d>1$, we substitute the asymptotic expansion (\ref{bssl}) with
$M=[{(d-1)}/{2}]+1$ for $J_{d/2-1}$ to write
\begin{equation}\label{hsm}
{\mathcal{H}}_{j}(\xi)=
{\mathcal{K}}_{j}(\xi)+
\frac{2^{(\frac{d}{2} +1)j}}{
|\xi|^{d}}\int_{\left\{r> 2^{-j}\right\}}
\Phi_{j}(r,\xi)
 {L_{\frac{d}{2}-1}(2^{j}r;
M)}r^{\frac{d}{2}}dr,
\end{equation}
where
\begin{equation*}
{\mathcal{K}}_{j}(\xi):=
\frac{2^{\frac{d+1}{2} j}}{|\xi|^{d}}
\sum_{\ell=0}^{M}
c_{\ell}(\tfrac{d}{2}-1){2^{-j\ell}}\int_{\left\{r> 2^{-j}\right\}}
r^{\frac{d-1}{2} -\ell}
\Phi_{j}(r,\xi)
\cos{(2^{j} r+\lambda_{\ell})}\,dr.
\end{equation*}
In the dimension $d=1$,
we obviously use (\ref{jd1}) to substitute for
$J_{-{1}/{2}}$ as opposed to the expansion
(\ref{bssl}). In that case, ${\mathcal{H}}_{j}(\xi)$ is simply the first term (i. e.,  the term that corresponds to $\ell=0$) of the latter sum that defines ${\mathcal{K}}_{j}(\xi)$. We are going to show that
\begin{equation}\label{hfrqhj}
\|{\mathcal{H}}_{j}\|_{
L^{p}(\mathbb{R}^{d})} \lesssim
{2^{d(1-\frac{1}{p})j+\delta(j)}},
\end{equation}
for all $1\leq p\leq \infty$. \par First we consider the sum ${\mathcal{K}}_{j}(\xi)$. Fix $0\leq \ell \leq M$. Recall that
\begin{equation*}
\Phi_{j}(r,\xi)=
  \psi\left(\frac{r}{2\pi |\xi|}\right)
E_{\alpha,\beta}\left( e^{\dot{\imath}\pi s} \left(\frac{2^{j}r}{2\pi  |\xi|}\right)^{\gamma}\right).
\end{equation*}
Since $\psi$ is smooth and has compact support and $r\mapsto E_{\alpha,\beta}(r^{\gamma})$ is smooth away from zero, one has
\begin{equation*}
\lim_{r\rightarrow +\infty} r^{a}|\partial^{m}_{r}\Phi_{j}(r,\xi)|=0,\qquad
\xi\neq 0,
\end{equation*}
for any real number $a$ and every positive integer $m$. Thus, integration by parts $m$ times we see that
\begin{multline}\label{hfrprts}
\int_{\left\{r> 2^{-j}\right\}}
r^{\frac{d-1}{2} -\ell}
\Phi_{j}(r,\xi)
\cos{(2^{j} r+\lambda_{\ell})}\,dr\\
=\sum_{k=1}^{m} (-1)^{k} {2^{-j k}}\sin{(1+\lambda_{\ell}+\pi k/2)} \left(\partial^{k-1}_{r}\left(
r^{\frac{d-1}{2} -\ell}
\Phi_{j}(r,\xi)\right)\right)_{r=2^{-j}}
 \\
+(-1)^{m}2^{-m j}
\int_{\left\{r> 2^{-j}\right\}}
\partial^{m}_{r}\left(
r^{\frac{d-1}{2} -\ell}
\Phi_{j}(r,\xi)\right)
\cos{(2^{j} r+\lambda_{\ell}+m\pi/2)}\,dr.
\end{multline}
We need to estimate the derivatives
on the right side of (\ref{hfrprts}).
As an application of formula (\ref{dme0}), one can show by induction that, for $r>0$, we have
\begin{equation*}
\partial_{r}^{m}
E_{\alpha,\beta}\left( e^{\dot{\imath}\pi s} \left(\frac{2^{j}r}{2\pi  |\xi|}\right)^{\gamma}\right)
=r^{-m} \sum_{k=1}^{m} a_{k}
e^{\dot{\imath}\pi k s}
\left(\frac{2^{j}r}{2\pi  |\xi|}\right)^{\gamma k}
\int_{C_{\omega}}
\frac{e^{\zeta^{1/\alpha}}{
\zeta^{(1-\beta)/\alpha}} d\zeta}{
\left(\zeta-e^{\dot{\imath}\pi s} \left(\frac{2^{j}r}{2\pi  |\xi|}\right)^{\gamma}\right)^{k+1}},
\end{equation*}
where $a_{k}$ are constants
that depend only on $\alpha$ and $\gamma$.
Using this we compute
\begin{multline*}
\partial^{m}_{r}\left(
r^{\frac{d-1}{2} -\ell}
\Phi_{j}(r,\xi)\right)
=\partial^{m}_{r}\left(r^{\frac{d-1}{2} -\ell} \psi\left(\frac{r}{2\pi |\xi|}\right)
E_{\alpha,\beta}\left( e^{\dot{\imath}\pi s} \left(\frac{2^{j}r}{2\pi  |\xi|}\right)^{\gamma}\right)
\right)\\
=\sum_{k_{1}+k_{2}+k_{3}=m}\sum_{k=k_{*}(k_{3})}^{k_{3}} c_{k_{1},k_{2},k_{3};m}a_{k} {r^{\frac{d-1}{2} -\ell-k_{1}-k_{3}}}{ |\xi|^{-k_{2}}}{\psi}^{(k_{2})}
\left(\frac{r}{2\pi |\xi|}\right) \\
e^{\dot{\imath}\pi k s}
\left(\frac{2^{j}r}{2\pi  |\xi|}\right)^{\gamma k} \int_{C_{\omega}}
\frac{e^{\zeta^{1/\alpha}}{
\zeta^{(1-\beta)/\alpha}}}{
\left(\zeta-e^{\dot{\imath}\pi s} \left(\frac{2^{j}r}{2\pi |\xi|}\right)^{\gamma}\right)^{k+1}} d\zeta,
\end{multline*}
with
$c_{k_{1},k_{2},k_{3};m}:=
\tfrac{m!}{{k_{1}!}{k_{2}!}{k_{3}!}}
(\tfrac{d-1}{2} -\ell)_{k_{1}}
(2\pi)^{-k_{2}}$, $a_{0}=(2 \pi \dot{\imath} \alpha)^{-1}$, and
\begin{equation*}
k_{*}(k):=
   \begin{cases}
    1, & k\geq 1; \\
      0, & k=0.
\end{cases}
\end{equation*}
Assume for the moment that
\begin{multline}\label{smt0}
\left(\frac{2^{j}r}{2\pi  |\xi|}\right)^{\gamma k} \left|{\psi}^{(k_{2})}
\left(\frac{r}{2\pi |\xi|}\right)\int_{C_{\omega}}
\frac{e^{\zeta^{1/\alpha}}{\zeta^{(1-\beta)/
\alpha}}}{
\left(\zeta-e^{\dot{\imath}\pi s} \left(\frac{2^{j}r}{2\pi |\xi|}\right)^{\gamma}\right)^{k+1}} d\zeta\right|\\ \lesssim 2^{\delta(j)}\left|{\psi}^{(k_{2})}
\left(\frac{r}{2\pi |\xi|}\right)\right|,
\end{multline}
for all $k\geq 0$, with an implicit constant independent of $r$ and $\xi$.
This would yield the estimate
\begin{equation}\label{stdr}
\left|\partial^{m}_{r}\left(
r^{\frac{d-1}{2} -\ell}
\Phi_{j}(r,\xi)\right)\right|
 \lesssim \,2^{\delta(j)}\sum_{k_{1}+k_{2}+k_{3}=m}\frac{r^{\frac{d-1}{2} -\ell-k_{1}-k_{3}}}{|\xi|^{k_{2}}}\left|{\psi}^{(k_{2})}
\left(\frac{r}{2\pi |\xi|}\right)\right|.
\end{equation}
Returning to the integral on the right side of (\ref{hfrprts}) with the estimate (\ref{stdr}), we find that
\begin{multline}\label{rtoxi}
\left|\int_{\left\{r> 2^{-j}\right\}}
\partial^{m}_{r}\left(
r^{\frac{d-1}{2} -\ell}
\Phi_{j}(r,\xi)\right)
\cos{(2^{j} r+\lambda_{\ell}+m\pi/2)}\,dr\right|\\
\lesssim \, {2^{\delta(j)}}\sum_{k_{1}+k_{2}+k_{3}=m}
\int_{\left\{r> 2^{-j}\right\}}
{r^{\frac{d-1}{2} -\ell-k_{1}-k_{3}}}{ |\xi|^{-k_{2}}}\left|{\psi}^{(k_{2})}
\left(\frac{r}{2\pi |\xi|}\right)\right| dr\\
\lesssim \, {2^{\delta(j)}}\sum_{k_{2}=0}^{m}
\int_{\left\{r> 2^{-j}\right\}}
{r^{\frac{d-1}{2} -\ell-m+k_{2}}}{ |\xi|^{-k_{2}}}\left|{\psi}^{(k_{2})}
\left(\frac{r}{2\pi |\xi|}\right)\right| dr.
\end{multline}
Using  (\ref{stdr}) once more, we see that the boundary terms in (\ref{hfrprts}) are bounded
by
\begin{multline}\label{stbndr}
\sum_{k=1}^{m} {2^{-k j}}\left|\left(\partial^{k-1}_{r}\left(
r^{\frac{d-1}{2} -\ell}
\Phi_{j}(r,\xi)\right)\right)_{r=2^{-j}}
 \right|\\
\lesssim
\,2^{\delta(j)}\,\sum_{k=1}^{m} {2^{-k j}}\sum_{k_{1}+k_{2}+k_{3}=k-1}{
2^{-j(\frac{d-1}{2} -\ell-k_{1}-k_{3})}}{ |\xi|^{-k_{2}}}\left|{\psi}^{(k_{2})}
\left(\frac{2^{-j}}{2\pi |\xi|}\right)\right|\\
\lesssim\,2^{\delta(j)-j(\frac{d-1}{2} -\ell+1)}\,\sum_{k=1}^{m} \sum_{k_{2}=0}^{k-1}{
2^{-j k_{2}}}{ |\xi|^{-k_{2}}}\left|{\psi}^{(k_{2})}
\left(\frac{2^{-j}}{2\pi |\xi|}\right)\right|.
\end{multline}
Note here that, for all $k_{2}\geq 0$, we have
\begin{equation}\label{rtoxi0}
\left.\begin{aligned}
\int_{\mathbb{R}^{d}}{
|\xi|^{-(d+k_{2})p}}\left|{\psi}^{(k_{2})}
\left(\frac{r}{2\pi |\xi|}\right)\right|^{p}
d\xi&\lesssim r^{d-(d+k_{2})p},\quad 1\leq p<\infty,\\
\sup_{\xi\in \mathbb{R}^{d}}{{
|\xi|^{-(d+k_{2})}}\left|{\psi}^{(k_{2})}
\left(\frac{r}{2\pi |\xi|}\right)\right|}
&\lesssim r^{-(d+k_{2})}.
\end{aligned}\right\}
\end{equation}
Now, multiply (\ref{rtoxi}) through by $|\xi|^{-d}$, then use the triangle inequality followed by Minkowski's integral inequality in the light of the estimates (\ref{rtoxi0}) to see that
\begin{multline}\label{nrm1}
\left\||\xi|^{-d}\int_{\left\{r> 2^{-j}\right\}}
\partial^{m}_{r}\left(
r^{\frac{d-1}{2} -\ell}
\Phi_{j}(r,\xi)\right)\cos{(2^{j} r+\lambda_{\ell}+m\pi/2)}\,dr\right \|_{L^{p}(\mathbb{R}^{d})}\\
\lesssim {2^{\delta(j)}}
\int_{\left\{r> 2^{-j}\right\}}
{r^{\frac{d}{p}-\frac{d+1}{2} -\ell-m}} dr\\
\approx {2^{\delta(j)-(\frac{d}{p}-\frac{d-1}{2} -\ell-m)j}},
\end{multline}
for all $1\leq p\leq \infty$, provided that $m>\frac{d}{p}-\frac{d-1}{2}$. Applying the same process to (\ref{stbndr}) shows that
\begin{equation}\label{nrm2}
\left\| |\xi|^{-d}
\sum_{k=1}^{m} {2^{-k j}}\left(\partial^{k-1}_{r}\left(
r^{\frac{d-1}{2} -\ell}
\Phi_{j}(r,\xi)\right)\right)_{r=2^{-j}}\right
\|_{L^{p}(\mathbb{R}^{d})}\,\lesssim\,
{2^{\delta(j)-(\frac{d}{p}-\frac{d-1}{2} -\ell)j}}.
\end{equation}
Combining the estimates (\ref{nrm1}) and (\ref{nrm2}), it follows from (\ref {hfrprts}) and the triangle inequality that
\begin{equation*}
\left\||\xi|^{-d}\int_{\left\{r> 2^{-j}\right\}}
r^{\frac{d-1}{2} -\ell}
\Phi_{j}(r,\xi)
\cos{(2^{j} r+\lambda_{\ell})}\,dr\right\|_{L^{p}(\mathbb{R}^{d})}\,\lesssim\,{2^{\delta(j)-(\frac{d}{p}-\frac{d-1}{2} -\ell)j}}.
\end{equation*}
By the triangle inequality, this estimate implies
\begin{equation}\label{hfrqhj0}
\|{\mathcal{K}}_{j}\|_{L^{p}(\mathbb{R}^{d})} \lesssim
{2^{d(1-\frac{1}{p})j+\delta(j)}},
\end{equation}
for all $1\leq p \leq \infty$.
To complete the proof of (\ref{hfrqhj0}), we need to verify (\ref{smt0}).
Fix $k\geq 0$ and $\xi \in \mathbb{R}^{d}\setminus\{0\}$ and let $j>0$. It follows from (\ref{nf1}) that, for all $r\in \operatorname{supp}{\psi}\left(\frac{r}{2\pi |\xi|}\right)$, we have
\begin{equation*}
\left(\frac{2^{j}r}{2\pi  |\xi|}\right)^{\gamma k}\sup_{\;\zeta\in C_{\omega}}{\left|\zeta-e^{\dot{\imath}\pi s} \left(\frac{2^{j}r}{2\pi |\xi|}\right)^{\gamma}\right|^{-(k+1)}}\lesssim
\left(\frac{2^{j}r}{2\pi  |\xi|}\right)^{-\gamma}\lesssim 2^{-\gamma j}.
\end{equation*}
It also follows from (\ref{nf1})
that, whenever $r\in \operatorname{supp}{\psi}\left(\frac{r}{2\pi |\xi|}\right)$, the left side of the latter inequality is bounded above by a constant independent of $j$ when $j<0$. The left side of (\ref{smt0}) is therefore majorized by a constant multiple of
\begin{equation*}
2^{\delta(j)} \left|{\psi}^{(k_{2})}
\left(\frac{r}{2\pi |\xi|}\right)\right|\,\int_{C_{\omega}}
{e^{|\zeta|^{1/\alpha}\cos{({\arg{\zeta}}/{\alpha
})}}{|\zeta|^{(1-\beta)/
\alpha}}} d\zeta.
\end{equation*}
This contour integral exists. Recall that ${\pi}/{2}<{\omega}/{\alpha}<\min\{\pi,{\pi}/
{\alpha}\}$.
Thus $\cos({\omega}/{\alpha})<0$ and the functions
$r\mapsto e^{r^{1/\alpha}\cos{({\pm \omega}/{\alpha
})}}{r^{(1-\beta)/
\alpha}}$ are  consequently in $L^{1}([1,\infty[)$. This proves the estimate (\ref{smt0}). An equivalent way to obtain (\ref{smt0}) is from the pointwise bound (see e.g. \cite{gorenflo2020mittag}, Section 4.4)
\begin{equation*}
|\partial_{r}^{k}
E_{\alpha,\beta}( e^{\dot{\imath}\pi s} r)|\lesssim
\frac{1}{1+r^{k+1}},\quad r>0,\,|s|>\alpha \pi/2.
\end{equation*}
 \par
It remains to estimate the last term on the right side of (\ref{hsm}). Since
\begin{equation*}
\left|{L_{\frac{d}{2}-1}(2^{j}r;
[\tfrac{d-1}{2}]+1 )}\right|\lesssim
(2^{j}r)^{-[\frac{d-1}{2} ]-\frac{5}{2}},
\end{equation*}
the remainder term in (\ref{hsm}) is bounded, modulo a multiplicative constant, by
\begin{multline}\label{rm1}
{2^{(d_{*}-1)j}}\int_{\left\{r> 2^{-j}\right\}}
|\xi|^{-d}\left|\Phi_{j}(r,\xi)\right|
r^{d_{*}-2}dr\\
\lesssim
{2^{(d_{*}-1)j+\delta(j)}}\int_{\left\{r> 2^{-j}\right\}}
|\xi|^{-d}\psi\left(\frac{r}{2\pi |\xi|}\right)
r^{d_{*}-2}dr
\end{multline}
by (\ref{stmml1}), where $d_{*}=\frac{d-1}{2}-[\frac{d-1}{2}]$. Obviously $0\leq d_{*}<1$. Using Minkowski's inequality, it follows from
(\ref{rtoxi0}) with $k_{2}=0$ that the
$L^{p}(\mathbb{R}^{d})$ norm of
the right side of (\ref{rm1}) is bounded by a constant multiple of
\begin{equation*}
{2^{(d_{*}-1)j+\delta(j)}}\int_{\left\{r> 2^{-j}\right\}}
r^{\frac{d}{p}-d+d_{*}-2}dr \,\approx\,
{2^{d(1-\frac{1}{p})j+\delta(j)}}.
\end{equation*}
Combining this with the estimate (\ref{hfrqhj0}) concludes the proof of (\ref{hfrqhj}). \par
Using the triangle inequality, we see from (\ref{plttrngl}) together with (\ref{nrmstl1}) and (\ref{hfrqhj}) that
\begin{equation}\label{pjstm0}
\|P_{j}\left(E_{\alpha,\beta}(
e^{\dot{\imath}\pi s}
 |\cdot|^{\gamma} )\right)^{\vee}
\|_{L^{p}(\mathbb{R}^{d})}\lesssim
{2^{d(1-\frac{1}{p})j+\delta(j)}}.
\end{equation}
Suppose that $\gamma>d(1-\frac{1}{p})>0$. Then the series
\begin{equation}\label{cnvrsr1}
\sum_{j\in \mathbb{Z}}
{{2^{d\left(1-\frac{1}{p}\right)j+\delta(j)}}}=
\sum_{j>0}
{{2^{\left(d\left(1-\frac{1}{p}\right)
-\gamma\right)j}}}+\sum_{j\leq0}
{{2^{d\left(1-\frac{1}{p}\right)j}}}
\end{equation}
is convergent. Finally, applying the Littlewood-Paley
inequality (\ref{plt}), the estimate (\ref{pjstm0})
together with the convergence of
the series (\ref{cnvrsr1}) yields Theorem \ref{lmlwp}.
\section{Applications}\label{apls}
Let $0<\alpha<2$, $\beta>0$, and $-1<\mu, \nu\leq 1$. Consider the following Cauchy problem
associated with a class of inhomogeneous space-time fractional equations:
 \begin{equation}\label{prblm1}
\left\{
 \begin{aligned}
e^{\mu \pi{\dot{\imath}}}\partial^{\alpha}_{t}u
 (t,x)&= e^{\nu \pi{\dot{\imath}}} (-\Delta)^{\beta/2}u(t,x)+
 F(t,x), \;\;(t,x)\in (0,\infty)\times {\mathbb{R}}^{d},\\
 u(0,x)&=f(x),\quad x\in{\mathbb{R}}^{d},
\end{aligned}
\right.
\end{equation}
where $\partial^{\alpha}_{t}$ denotes
the Caputo derivative of order $\alpha$ defined by
\begin{equation*}
\partial^{\alpha}_{t}U(t):=
\frac{1}{\Gamma{(1-\alpha)}}
\int_{0}^{t}
\frac{U^{\prime}(\tau)}{(t-\tau)^{\alpha}}d\tau,
\end{equation*}
and $(-\Delta)^{\beta/2}$
is the Fourier multiplier with symbol $|\xi|^{\beta}$ well-known as the fractional Laplacian. \par
When $\mu=0$ and $\nu=1$, the problem (\ref{prblm1})
coincides with the space-time fractional heat equation studied in   \cite{chen2012space,kemppainen2017representation,li2022applicationi,
li2022applicationii,mainardi2001fundamental}.
If $\mu=\alpha/2$ and $\nu=0$ then (\ref{prblm1}) reduces to the space-time fractional Schr\"{o}dinger equation investigated in \cite{bayin2013time,dong2008space,grande2019space,
lee2020strichartz,su2020holder,wang2007generalized}. Moreover, if $\mu=1/2$ and $\nu=0$, then (\ref{prblm1}) simplifies to the
space-time fractional Schr\"{o}dinger equation studied in \cite{su2019local,su2021dispersive}. For the physical interpretation and significance of these
equations we refer e.g. to \cite{laskin2000fractional,meerschaert2002stochastic}.\par
For regular enough $f$ and $F$, an explicit formula for the solution of (\ref{prblm1}) takes the form:
\begin{equation}\label{sol11}
u(t,x)=v(t,x)+w(t,x),
\end{equation}
where
\begin{eqnarray*}
v(t,x)&:=&V_{\alpha,\mu,\nu}(t,\cdot)\ast
f(x),\\
w(t,x)&:=&e^{-\mu \pi{\dot{\imath}}}\int_{0}^{t}
(t-s)^{\alpha-1}W_{\alpha,\mu,\nu}(t-s,\cdot)\ast
{F(s,\cdot)}(x)ds,
\end{eqnarray*}
with
\begin{eqnarray*}
V_{\alpha,\mu,\nu}(t,x) &:=& \int_{\mathbb{R}^{d}}
e^{2\pi \dot{\imath}x\cdot \xi}
E_{\alpha}\left(e^{(\nu-\mu) \pi {\dot{\imath}}}t^{\alpha} |\xi|^{\beta}
\right)d\xi,\\
W_{\alpha,\mu,\nu}(t,x) &:=& \int_{\mathbb{R}^{d}}
e^{2\pi \dot{\imath}x\cdot \xi}
E_{\alpha,\alpha}\left(e^{(\nu-\mu) \pi {\dot{\imath}}}t^{\alpha} |\xi|^{\beta}
\right)d\xi.
\end{eqnarray*}
The solution (\ref{sol11}) is obtained rigorously in Section \ref{plctfrm} below. It is evidently consistent with Duhamel's principle since $v$ solves the homogeneous problem
\begin{equation}\label{hmpr1}
\left\{
 \begin{aligned}
e^{\mu \pi{\dot{\imath}}}\partial^{\alpha}_{t}v
 (t,x)&= e^{\nu \pi{\dot{\imath}}} (-\Delta)^{\beta/2}v(t,x), \;\;(t,x)\in (0,\infty)\times {\mathbb{R}}^{d},\\
 v(0,x)&=f(x),\quad x\in{\mathbb{R}}^{d},
\end{aligned}
\right.
\end{equation}
and $w$ solves the inhomogeneous problem
\begin{equation}\label{hmpr2}
\left\{
 \begin{aligned}
e^{\mu \pi{\dot{\imath}}}\partial^{\alpha}_{t}w
 (t,x)&= e^{\nu \pi{\dot{\imath}}} (-\Delta)^{\beta/2}w(t,x)+F(t,x), \;\;(t,x)\in (0,\infty)\times {\mathbb{R}}^{d},\\
 w(0,x)&=0,\quad x\in{\mathbb{R}}^{d}.
\end{aligned}
\right.
\end{equation}
Our first goal is to find the largest range of values of the exponents $p$, $q$ such that the dispersive estimate (\ref{hmdsrs1})
of Theorem \ref{thds} below holds for all $\alpha \in (0,2)$, $\beta>0$ and $\mu,\nu\in (-1,1]$ such that $|\mu-\nu|>\alpha/2$. Before the recent work of Kemppainen et al. \cite{kemppainen2017representation} and
Su et al. \cite{su2019local}, this estimate was known to hold only in the high regularity regime $\beta > d/2$, for certain values of $p$ and $q$, in the cases of the heat equation \cite{kemppainen2017representation,li2022applicationi}
and the Schr\"{o}dinger equation \cite{su2021dispersive}. An estimate similar to (\ref{hmdsrs1}) is deduced for the convolution $W_{\alpha,\mu,\nu}(t-s,\cdot)\ast
{F(s,\cdot)}$, $0<s<t$, and is used to prove a
decay estimate for the $L^{p}(\mathbb{R}^d)$ norm
of the solution of the inhomogeneous problem (\ref{hmpr2}) in Theorem \ref{thnhm}.
One can adapt the estimates proved here to obtain
local well-posedness results analogous to those proved in \cite{su2019local} for the semilinear problem
 \begin{equation*}
\left\{
 \begin{aligned}
e^{\mu \pi{\dot{\imath}}}\partial^{\alpha}_{t}u
 (t,x)&= e^{\nu \pi{\dot{\imath}}} (-\Delta)^{\beta/2}u(t,x)\pm
|u(t,x)|^{\sigma}\,u(t,x), \;(t,x)\in (0,\infty)\times {\mathbb{R}}^{d},\\
 u(0,x)&=f(x),\quad x\in{\mathbb{R}}^{d},\;\;
\sigma>0.
\end{aligned}
\right.
\end{equation*}
\subsection{An explicit formula for the solution}
\label{plctfrm}
\par Let $\widehat{\upsilon}$ and $\widetilde{\upsilon}$ denote the Fourier and Laplace transforms of a given nice enough function $\upsilon$, respectively. Assume that $f$ and $F(t,\cdot)$ are Schwartz functions on ${\mathbb{R}^{d}}$ and let $F(\cdot,\xi)$ have a well-defined Laplace transform. Taking the Fourier transform of (\ref{prblm1}) in the spacial variable, it transforms into the fractional differential equation
\begin{equation}\label{frlp}
 e^{\mu \pi{\dot{\imath}}}
 \partial^{\alpha}_{t}\widehat{u}(t,\xi)-
  e^{\nu \pi{\dot{\imath}}}
|\xi|^{\beta}
\widehat{u}(t,\xi)=\widehat{F}(t,\xi),
\end{equation}
with the initial condition $\widehat{u}(0,\cdot)=
\widehat{f}$. Taking the Laplace transform
of both sides of (\ref{frlp}) in $t$
and using the fact that
$\widetilde{\partial^{\alpha}_{t}U}(s)=
s^{\alpha}\widetilde{U}(s)-s^{\alpha-1}U(0)$
yield
\begin{equation*}
\left(s^{\alpha}-e^{(\nu-\mu)\pi {\dot{\imath}}}|\xi|^{\beta}\right)
\widetilde{\widehat{u}}(s,\xi)=
{s^{\alpha-1}}\widehat{f}(\xi)
+e^{-\mu \pi{\dot{\imath}}}{\widetilde{\widehat{F}}(s,\xi)}.
\end{equation*}
Solving this algebraic equation
for $\widetilde{\widehat{u}}(s,\xi)$ then inverting the Laplace transform we obtain
\begin{multline}\label{sol1}
\widehat{u}(t,\xi)=
E_{\alpha}\left(e^{(\nu-\mu) \pi {\dot{\imath}}}t^{\alpha} |\xi|^{\beta}
\right)\widehat{f}(\xi)\\
+ e^{-\mu \pi{\dot{\imath}}}\int_{0}^{t}(t-s)^{\alpha-1}{\widehat{F}(s,\xi)}
E_{\alpha,\alpha}\left(e^{(\nu-\mu) \pi {\dot{\imath}}}(t-s)^{\alpha} |\xi|^{\beta}
\right)ds.
\end{multline}
One can verify  that (\ref{sol1}) solves (\ref{frlp}), mainly, by computing the Laplace transform of the series
that represent $E_{\alpha}(t^{\alpha})$ and $t^{\alpha-1}E_{\alpha,\alpha}(t^{\alpha})$
term by term. For more details on solving Cauchy problems with Caputo derivatives, see
Theorem 7.2 in \cite{diethelm2010analysis},
Sections 7.2.1 and 7.3.2 in \cite{gorenflo2020mittag}, and Appendix A in \cite{grande2019space}.
\subsection{Dispersive estimates}
By radial symmetry and a scaling argument, the kernel of the solution of the homogeneous problem
(\ref{hmpr1}) takes the form
\begin{equation*}
V_{\alpha,\mu,\nu}(t,x) =
t^{-\frac{\alpha}{\beta}d}(
E_{\alpha}(e^{(\nu-\mu) \pi {\dot{\imath}}} |\cdot|^{\beta}))^{\widehat{}}
\,(t^{-\frac{\alpha}{\beta}}x).
\end{equation*}
Let $1\leq p,q\leq \infty$, and let $1-\frac{1}{r}=\frac{1}{p}-\frac{1}{q}$.
Invoking Young's convolution inequality, we obtain \begin{equation*}
\|v(t,\cdot)\|_{{L^{q}(\mathbb{R}^{d})}}\leq
{t^{-\frac{\alpha}{\beta}d(1-\frac{1}{r})}}
\|(
E_{\alpha}(e^{(\nu-\mu) \pi {\dot{\imath}}} |\cdot|^{\beta}))^{\widehat{}}
\,\|_{L^{r}(\mathbb{R}^{d})}
\|f\|_{L^{p}(\mathbb{R}^{d})},
\end{equation*}
whence, by Theorem \ref{lmlwp}, we have\smallskip
\begin{theorem}\label{thds}
Assume that $\alpha\in (0,2)$ and
$\,|\mu-\nu|>\alpha/2$. Then the dispersive estimate
\begin{equation}\label{hmdsrs1}
\|v(t,\cdot)\|_{{L^{q}(\mathbb{R}^{d})}}\lesssim
{t^{-\frac{\alpha}{\beta}d(\frac{1}{p}-\frac{1}{q})}}
\|f\|_{L^{p}(\mathbb{R}^{d})},
\end{equation}
holds true for all $1\leq p< q\leq \infty$ such that
$\frac{1}{p}-\frac{1}{q}<\frac{\beta}{d}$.
\end{theorem}
\subsection{Estimates for the inhomogeneous problem}
Let $w$ be a solution of the inhomogeneous problem (\ref{hmpr2}). Using Minkowski's inequality, we see that
\begin{equation}\label{wtxlq}
\|w(t,\cdot)\|_{{L^{q}(\mathbb{R}^{d})}}
\leq \int_{0}^{t}
(t-s)^{\alpha-1}\|W_{\alpha,\mu,\nu}(t-s,\cdot)\ast
{F(s,\cdot)}\|_{{L^{q}(\mathbb{R}^{d})}} ds.
\end{equation}
Analogously to the estimate (\ref{hmdsrs1}), if
$\alpha \in (0,2)$ and $\,|\mu-\nu|> \alpha/2$ then
\begin{equation*}
\|W_{\alpha,\mu,\nu}(t-s,\cdot)\ast
{F(s,\cdot)}\|_{{L^{q}(\mathbb{R}^{d})}}\lesssim
{(t-s)^{-\frac{\alpha}{\beta}d(\frac{1}{p}-\frac{1}{q})}}
\|{F(s,\cdot)}\|_{{L^{p}(\mathbb{R}^{d})}},
\end{equation*}
for all $1\leq p< q\leq \infty$ such that
$\frac{1}{p}-\frac{1}{q}<\frac{\beta}{d}$.
Using this estimate in (\ref{wtxlq})
then applying H\"{o}lder's inequality, we deduce\smallskip
\begin{theorem}\label{thnhm}
Suppose that $\alpha\in (0,2)$, ${1}/{\alpha}<r\leq \infty$, and that $1\leq p<q\leq \infty$ are such that
$\frac{1}{p}-\frac{1}{q}<\frac{\beta}{d}
(1-\frac{1}{\alpha}\frac{1}{r})$. Then
\begin{equation*}
\|w(t,\cdot)\|_{{L^{q}(\mathbb{R}^{d})}}
\lesssim t^{\frac{\alpha}{\beta}d\left(\frac{\beta}{d}
-\left(\frac{1}{p}-\frac{1}{q}\right)\right)-\frac{1}{r}}
\|F\|_{L^{r}([0,\infty);L^{p}(\mathbb{R}^{d}))}.
\end{equation*}
\end{theorem}
\bigskip
\subsection*{Acknowledgements}
The author would like to thank the reviewers for their useful comments and the editors for their support.
\bibliography{24027FIN}

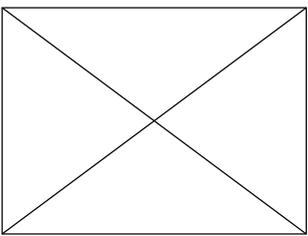
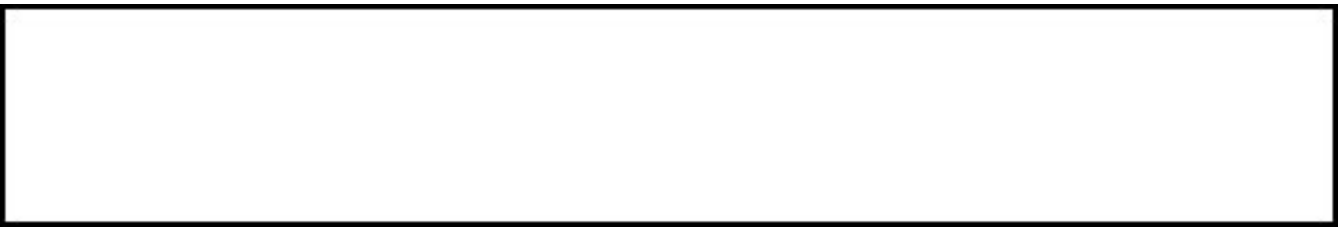
\end{document}